\documentclass[journal, onecolumn]{IEEEtran}
\ifCLASSINFOpdf
\else
\fi
\usepackage{amsopn}
\RequirePackage{amsmath,amssymb,mathrsfs,bbm,amsthm}
\RequirePackage{mathtools}
\usepackage{lipsum}
\usepackage{amsfonts}
\usepackage{graphicx}
\usepackage{epstopdf}
\usepackage[colorlinks]{hyperref}
\usepackage[export]{adjustbox}
\usepackage{cite}

\newcommand{\mc}{\mathcal}
\newcommand{\ms}{\mathscr}
\newcommand{\tv }[1]{ d_{TV}\big( #1 \big)  }
\newcommand{\pois}[1]{ \text{Pois}(#1)}
\newcommand{\knn}[2]{ \rho_{i,#1,#2}  }
\newcommand{\ind}{\mathbbm{1} }
\DeclarePairedDelimiter{\ceil}{\lceil}{\rceil}
\DeclarePairedDelimiter{\floor}{\lfloor}{\rfloor}

\newtheorem{theorem}{Theorem}
\newtheorem{lemma}{Lemma}

\begin{document}
\title{Consistent Entropy Estimation for Stationary Time Series \thanks{Submitted to the editors \today.}}
%\funding{Supported by grants \#1045153 and  \#1546130 from the National Science Foundation}}}

\author{
\thanks{This work was supported by grants \#1045153 and  \#1546130 from the National Science Foundation. Portions of this paper were presented as a poster at the 2019 North American School on Information Theory.}
Alexander L. Young\thanks{A.L. Young is with the Department of Statistics, Harvard University, Cambridge, MA 02138 USA
  (email: \href{mailto:alexander_young@fas.harvard.edu}{alexander\_young@fas.harvard.edu}).} and 
\and David B. Dunson \thanks{D.B. Dunson is with the Department of Statistical Sciences, Duke University, Durham, NC 27705 USA} }


\markboth{IEEE TRANSACTIONS ON INFORMATION THEORY: submitted \today}{}%



\maketitle

\begin{abstract}
Entropy estimation, due in part to its connection with mutual information, has seen considerable use in the study of time series data including causality detection and information flow.  In many cases, the entropy is estimated using $k$-nearest neighbor (Kozachenko-Leonenko) based methods. However, analytic results on this estimator are limited to independent data.  In the article, we show rigorous bounds on the rate of decay of the bias in the number of samples, $N$, assuming they are drawn from a stationary process which satisfies a suitable mixing condition.   Numerical examples are presented which demonstrate the efficiency of the estimator when applied to a Markov process with stationary Gaussian density. These results support the asymptotic rates derived in the theoretical work.
\end{abstract}

\begin{IEEEkeywords}
Entropy estimation,
information theory, 
nearest neighbors, 
nonparametric estimation, 
stationary processes
\end{IEEEkeywords}

\section{Introduction}
\label{SEC:Intro}
For a random variable $X \in \mathbb{R}^d$ with density $f:\mathbb{R}^d \to [0,\infty)$, the (differential) entropy of $X$ is
\begin{equation}
H(X) = H(f) = - E[\log f(X)] = - \int f(x) \log f(x) dx.
\end{equation}
Here we have adopted the convention $0\log 0 = 0$ in the case that $f$ is not supported on all of $\mathbb{R}^d.$

Entropy is the main quantity of interest in information theory and has numerous applications in statistics \cite{berrett2019} making its estimation from data a desirable goal. One important case motivating this article is the estimation of mutual information, a measure of dependence capturing any relationship between random variables.  Given random variables $(X,Y)$ the mutual information is defined in terms of the joint and marginal entropies of $X$ and $Y$ as
\begin{equation}
\label{EQ:Def_MI}
I(X,Y) = H(X) +H(Y) - H(X,Y).
\end{equation}
This quantity is zero when $X$ and $Y$ are independent and positive otherwise.  Mutual information has been a topic of considerable interest in causality detection and information flow in bivariate time series \cite{PALUS199664, dionisio_mutual_2004, frenzel_partial_2007, Keller:2015:EMI:2791347.2791348, liu_mutual_2016, lungarella_methods_2007, palus_directionality_2007, palus_synchronization_2001, doi:10.1142/S0218127409025298,pham_mutual_2002, schreiber_measuring_2000, vejmelka_inferring_2008}.  Thus, methods for entropy estimation in dependent data are of clear importance.

Unfortunately, given samples $X_0,\dots,X_N$ drawn from $f$, one cannot typically construct an estimator for $H(f)$ of the form 
\begin{equation}
\label{EQ:Naive_Est}
\frac{-1}{N+1}\sum_{i=0}^N \log f(X_i )
\end{equation}
as $f$ is typically unknown.  A myriad of strategies have been proposed including kernel-density and B-spline based estimators of $H(X)$ \cite{walters-williams_estimation_2009}. Our focus is on the popular and well studied Kozachenko-Leonenko estimator of $H(X)$ which uses the nearest neighbors of $X_i$ to estimate $f(X_i)$ nonparametrically \cite{berrett2019, delattre_kozachenkoleonenko_2017, kozachenko_statistical_1987,  ent_diff}.  

We begin with a brief review of this estimator.  Suppose we have $N+1$ samples of $\mathbb{R}^d$-valued random variables $X_0,\dots,X_N$ drawn from $f$.  For $i\in \{0,\dots,N\}$, let $\mc N_i = \{0,\dots,N\}\setminus \{i\}.$ We then define $\knn{k}{\mc N_i}$ to be the distance from $X_i$ to its $k$th nearest neighbor in the set $\{X_j\}_{j \in \mc N_i}$ under a metric $\| \cdot \|$ defined on $\mathbb{R}^d$. Note there are $k$ points in the set $\{X_j\}_{j\in \mc N_i}$ contained in the closed ball centered at $X_i$ with radius $\knn{k}{\mc N_i}.$ Therefore, assuming regularity of $f$ and that $\knn{k}{\mc N_i}$ is small we may rely on the estimate \cite{delattre_kozachenkoleonenko_2017} 
$$ f(X_i) \knn{k}{\mc N_i}^d \nu_d \approx \frac{k}{N}, $$ where $\nu_d$ is the volume of the unit ball in $\mathbb{R}^d$ under metric $\| \cdot \|.$  This leads to the approximation 
$$ - \log f(X_i) \approx \log \frac{N \knn{k}{\mc N_i}^d \nu_d}{e^{\Psi(k)}}.$$
For small $k$, a bias correction replaces $k$ with $e^{\Psi(k)}$ where $\Psi(k)$ is the digamma function \cite{delattre_kozachenkoleonenko_2017}.  The Kozachenko-Leonenko or $k$-nearest neighbor entropy estimator of $H(f)$ is then
\begin{equation}
\hat{H}_{N+1}=  \frac{1}{N+1} \sum_{i=0}^N \log Y_i , \qquad Y_i=  \frac{N \knn{k}{\mc N_i}^d \nu_d}{e^{\Psi(k)}}.
\label{EQ:KLKNN_Estimator}
\end{equation}

Bias, variance, and central limit type theorems in the $k=1$ case were presented in \cite{delattre_kozachenkoleonenko_2017}. The $k>1$ case was considered in \cite{berrett2019} including a novel estimator using weighted averages of \eqref{EQ:KLKNN_Estimator} to improve efficiency in higher dimensions. Additionally, \cite{kraskov_estimating_2004} presented a bias-corrected version of \eqref{EQ:KLKNN_Estimator} specifically for the estimation of mutual information.  Theoretical guarantees for that mutual information estimator and additional analysis of the $k$-nearest neighbor entropy estimation may be found in \cite{demyst}. However, in each of these cases, the data were assumed to be independent and identically distributed (\emph{iid}).  

The primary focus of article, motivated by the performance of \eqref{EQ:KLKNN_Estimator} for \emph{iid} data and the desire for entropy estimation in time series, is to assess the bias of \eqref{EQ:KLKNN_Estimator} when applied to a stationary process $X_0,\dots,X_N$ with invariant density $f$.  Notably, we demonstrate that the Poisson approximation used to study \eqref{EQ:KLKNN_Estimator} in the \emph{iid} cases holds when $X_0,\dots, X_N$ satisfies a $\psi$-type mixing condition \cite{bradley2005} in addition to a H\"older continuity requirement on $f$ \cite{berrett2019, delattre_kozachenkoleonenko_2017}. We also present tail bounds on the distribution of $\knn{k}{\mc N_i}$ using a data-thinning approach in conjunction with asymptotics on order statistics. These results require bounds on the moments of $f$.

The remainder of this article is organized as follows: in Section \ref{SEC:Theory}, we introduce the relevant notation and present and discuss our primary assumptions.  In Section \ref{SEC:Results}, we introduce the main results and important lemmas.  Finally, in Section \ref{SEC:Example}, we compare the theoretical rates at which the bias decays with results from simulations of an autoregressive process which has an invariant Gaussian density. The proofs of all Lemmas are contained in an Appendix.

\section{Theory}
\label{SEC:Theory}
\subsection{Mathematical Formulation, Definitions, and Notation}

Let $\mathbb{Z}$ denote the integers, $\mathbb{R}^d$ denote $d$-dimensional Euclidean space, and $\mc B(\mathbb{R}^d)$ denote the Borel $\sigma$-algebra on $\mathbb{R}^d$.  Given $X\in\mathbb{R}^d$ we let $X^{(j)}$ denote the $j$th coordinate of $X$. We let $B(x,r) = \{ y\in \mathbb{R}^d: \| y-x\| < r\}$ denote the ball of radius $r$ centered at $x$ and $\ms L(W)$ denote the law of random variable $W$. Expectation is denoted by $E[\cdot]$ and $\ind_A(X)$ is the indicator for $X\in A$. Given two probability measures $\mu$ and $\lambda$ defined on the same probability space $(\Omega,\mc F)$, we define the total variation distance between $\mu$ and $\lambda$ to be 
\begin{equation}
\tv{\mu,\lambda} = \sup_{A \in \mc F} |\mu(A) - \lambda(A) |
\label{EQ:TV_def}
\end{equation}
We use $\pois{\lambda}$ to denote the Poisson distribution with rate $\lambda>0$ and $\Gamma(a,b)$ to denote the Gamma distribution with shape $a>0$ and rate $b>0$.

Given a multi-index $\gamma = (\gamma_1,\dots,\gamma_d)$ we let $|\gamma| = \sum_{i=1}^d \gamma_i$, $\gamma ! = \gamma_1! \dots \gamma_d!$, $X^\gamma = {X^{(1)}}^{\gamma_1}\dots{X^{(d)}}^{\gamma_d}$, and $D^\gamma f = \partial_{x_1}^{\gamma_1}\dots\partial_{x_d}^{\gamma_d} f.$ For a subset $U\subset \mathbb{R}^d$, we let $C^{\ell,\alpha}(\overline{U})$ be the H\"older space consisting of all functions from $U$ to $\mathbb{R}$ which are $\ell$-times continuously differentiable and have $\ell$th derivatives which are H\"older continuous with exponent $\alpha.$   A function $f:\mathbb{R}^d \to \mathbb{R}$ is in $C^{\ell,\alpha}(\overline{U})$ if the norm 
\begin{equation}
\|f\|_{C^{\ell,\alpha}(\overline{U})} = \sum_{|\gamma| \le \ell} \|D^\gamma f\|_{C(\overline{U})} + \sum_{|\gamma| = \ell} \|D^\gamma\|_{C^{0,\alpha}(\overline{U})}
\end{equation}
if finite \cite{evans} where 
$$\| f\|_{C(\overline{U})} = \sup_{x\in U} |f(x)|, \quad\|f\|_{C^{(0,\alpha)}(\overline{U})} = \sup_{x\ne y \in U} \frac{|f(x) - f(y)|}{\|x-y\|^\alpha}.$$  Hereafter, we let $C_f = \|f\|_{C^{\ell,\alpha}(\overline{U})}.$

Let $P$ be a probability measure on the infinite product space $\Omega = \prod_{i \in \mathbb{Z}} \mathbb{R}^d$ equipped with the product $\sigma$-algebra $\mc F = \otimes_{i \in \mathbb{Z}}  \mc B(\mathbb{R}^d).$  Let ${\bf X} = \{X_i; i \in \mathbb{Z}\}$ be a discrete-time stochastic process defined on $\Omega$ by the coordinate projections $X(i,\cdot) = X_i.$  For any subset $\mc I \subset \mathbb{Z}$, we let $\mc F_{\mc I} = \sigma(X_i; i\in \mc I) \subset \mc F$ denote the sub $\sigma$-algebra generated by $\{X_i\}_{i\in \mc I}.$  

Additionally, we assume that $\{X_i; \, i \in \mathbb{Z}\}$ is a stationary process so that for any $s,\tau \in \mathbb{Z}$ with $\tau >0$, $s\ge 1$ and $\{t_1,\dots,t_s\} \subset \mathbb{Z}$ 
\begin{equation}
\label{EQ:Stationary}
\ms L \big( (X_{t_1},\dots,X_{t_s})\big) = \ms L \big( (X_{t_1+\tau},\dots,X_{t_s+\tau})\big).
\end{equation}
This implies $\ms L(X_i)$ is equivalent for all $i \in \mathbb{Z}$ which we assume is absolutely continuous with respect to the Lebesgue measure on $\mathbb{R}^d$ with density $f.$ Thus, $\{X_i;\, i\in\mathbb{Z}\}$ are dependent, identically distributed random variables with common marginal distribution 
\begin{equation}
\label{EQ:Marginal}
P(X_i \in A ) = \int_A f_X(x) dx, \quad \forall A \in \mc B(\mathbb{R}^d).
\end{equation}

Of course, the dependency in our stochastic process should be expected to play an important role in the consistency of \eqref{EQ:KLKNN_Estimator}.  For technical reasons, we assume that $X_i\ne X_j$, $P$-almost surely for all $i\ne j$ to ensure that the $\knn{k}{\mc N_i} >0$.

\subsection{Assumptions}
There are a few very important assumptions necessary to demonstrate the consistency of \eqref{EQ:KLKNN_Estimator} in our setting. Namely,
\begin{itemize}
\item[A1] (Regularity) There exists $\alpha \in (0,1]$ such that the stationary density $f:\mathbb{R}^d \to [0,\infty)$ is in $C^{2,\alpha}(\overline{supp(f)})$
where $supp(f) = \{x: f(x) >0\}$ is the support of $f$. 
\item[A2] (Moments) There exists $r >0$ such that 
\begin{equation}
\int \| x\|^{d+r} f(x)dx < \infty.
\label{EQ:Moment}
\end{equation}

\item[A3] (Mixing) There exists a function $\psi : \mathbb{R} \to [0,\infty)$ and positive constants $K$ and $\epsilon$ satisfying the bound
\begin{equation}
\label{EQ:Mixing_bound}
\psi(z) \le \frac{K}{1+|z|^{1+\epsilon}}
\end{equation} 
such that for $i \notin \mc J \subset \mathbb{Z}$ and any events $A \in \mc B(\mathbb{R}^d)$ and $B \in  \mc F_{\mc J}$
\begin{equation}
\label{EQ:Mixing}
\frac{\big| P(X_i \in A \mid B ) - P(X_i \in A)\big|}{P(X_i \in A)} \le \psi\big( d(i,\mc J)\big) 
\end{equation}
such that $P(X_i \in A)  > 0$. Here $d(i,\mc J) = \min_{j \in \mc J} |i-j|.$ 

\end{itemize}

Both (A1) and (A2) are similar to the assumptions necessary in the \emph{iid} setting \cite{berrett2019, delattre_kozachenkoleonenko_2017}, and they have important implications.  First, under (A1) it follows that $f$ admits the following expansion,
\begin{equation}
f(y) = f(x) + \sum_{|\gamma| = 1}^2 \frac{1}{\gamma!} D^\gamma f(x) (y-x)^\gamma + C(x) \|y-x\|^{\alpha+2}
\label{EQ:Taylor_expansion}
\end{equation}
where $f(x)$, $C(x)$, and $D^\gamma f(x)$ are bounded above by $ C_f.$ Importantly, integrating this expansion over $B(x,r)$ provides the bound
\begin{equation}
\label{EQ:Holder_bound_on_integral}
\bigg|   \int_{B(x,r)} f(y)dy - f(x)\nu_d r^d  \bigg| \le \frac{d^2\nu_d C_f}{2} r^{2+d} + C_f\nu_d r^{2+d+\alpha}
\end{equation}
as the terms $\int_{B(x,r)} (y-x)^\gamma dy = 0 $ whenever $|\gamma|$ is odd \cite{delattre_kozachenkoleonenko_2017}. 
Secondly, under (A2) 
\begin{equation}
\int f^{1-\theta}(x) dx <\infty, \quad \forall \theta \in \bigg(0, \frac{d+r}{2d+r}\bigg).
\end{equation} 
Briefly, we may write 
$$
f^{1-\theta}(x) = \big( f(x) [1+\|x\|]^{(d+r)}\big)^{1-\theta} \frac{1}{[1+\|x\|]^{(d+r)(1-\theta)}}
$$
Taking $p = 1/(1-\theta)$ and $q=1/\theta$, it follows from H\"older's inequality that
\begin{align*}
\label{EQ:Moment_and_Holder}
\int f^{1-\theta}(x) dx \le &   \bigg(\int f(x)[1+\| x \| ]^{d+r}dx \bigg)^{1-\theta} \\
&  \times \bigg( \int \frac{1}{[1+\|x\| ]^{(1-\theta)(d+r)/\theta} } dx \bigg)^\theta
\end{align*}  
The right side of the above expression is finite since $\frac{(1-\theta)(r+d)}{\theta} > d$ and $\int [1+\|x\|]^{d+r} f(x) dx <\infty$. 

\subsubsection{Consistency in the Independent, Identically Distributed Case}

It will be useful to understand how the asymptotic analysis is performed in the \emph{iid} case prior to proving our main result for time series data. Please see \cite{berrett2019} for a more thorough exposition. Consider the distribution of $Y_i$ conditional on $X_i$ where $f(X_i)>0$.  For $r \in (0,\infty)$, consider the sequence of distribution functions
\begin{equation}
\begin{split}
F_{N,X_i}(r) &= P(Y_i \le r \mid X_i)  = P\big(\knn{k}{\mc N_i} \le r_{N} \mid X_i \big),\\ r_{N} &= \bigg(\frac{re^{\Psi(k)}}{N\nu_d} \bigg)^{1/d}.
\end{split}
\label{EQ:Estimator_CDF}
\end{equation}
Note that \eqref{EQ:Estimator_CDF} is equivalent to the requirement that at least $k$ observations from the set $\{X_j\}_{j\in \mc N_i}$ are within $r_{N,u}$ of $X_i$.  Assuming these random variables are independent, it follows that 
\begin{equation}
F_{N,X_i}(r) = \sum_{j=k}^N \binom{N}{j} p_{r,X_i,N}^j(1-p_{r,X_i,N})^{N-j}
\label{EQ:CDF_Independent_Binomial_Exact}
\end{equation}
where $$p_{r,X_i,N} = \int_{B(X_i, r_{N})} f(y)dy.$$  For fixed $r$, under (A1)
\begin{equation}
\begin{split}
p_{r,X_i,N} &\approx f(X_i) r_{N}^d \nu_d +O(r_{N}^{(2+d)/d})\\& = f(X_i) \frac{re^{\Psi(k)}}{N} + O\big(N^{-(d+2)/d}\big)
\end{split}
\end{equation}
as $N\to\infty$. 

Importantly, $p_{u,X_i}$ is $O(N^{-1})$ at leading order, motivating the use of a Poisson approximation in   \eqref{EQ:CDF_Independent_Binomial_Exact} as $N\to \infty$. Namely
\begin{equation}
F_{X_i}(r) = \lim_{N\to \infty} F_{n,X_i}(r) = \sum_{j=k}^\infty e^{-\lambda_{r,X_i}} \frac{\lambda_{r,X_i}^j}{j!} 
\label{EQ:CDF_Independent_Poisson_Approx}
\end{equation}
where $\lambda_{r,k,X_i} = f(X_i)re^{\Psi(k)}.$ Formally,
\begin{align*}
&E[\log Y_i \mid X_i]  = \int_0^\infty \log u \, dF_{N,x}(r) \\ 
&\xrightarrow{N\to\infty}  \int_0^\infty \log r\, dF_{X_i}(r) = \int_0^\infty \log r \frac{\lambda_{r,k,X_i}^{k-1}}{(k-1)!}e^{-\lambda_{r,k,X_i}}dr \\
\end{align*}
The final expression in the equation above is $E \log T_{k,X_i}$ where $T_{k,X_i} \sim \Gamma(k,e^{\Psi(k)}f(X_i))$ which is $-\log f(X_i).$
Thus,
\begin{equation}
E[\log Y_i] = E\big[E[\log Y_i \mid X_i] \big] \xrightarrow{N\to\infty} -E[\log f(X_i)] = H(f)
\end{equation}
suggesting that \eqref{EQ:KLKNN_Estimator} is a consistent estimator of $H(f)$. 

The fundamental challenge to validating the use of \eqref{EQ:KLKNN_Estimator} in time series is demonstrating the validity of the Poisson approximation in \eqref{EQ:CDF_Independent_Poisson_Approx} for dependent data.  There is considerable literature on this topic \cite{erhardsson_compound_1999, leadbetter_extremal_1988, leadbetter_extremes_1983, arratia1989, aldous_probability_1989}.  We will make use of the Stein-Chen method of \cite{arratia1989} in our analysis, which relies on bounding the first and second moments of the number of times the process $\{X_j\}_{j\in \mc N_i}$ visits a shrinking neighborhood of $X_i.$  

\section{Main results: efficiency of \eqref{EQ:KLKNN_Estimator} in the time series setting}
\label{SEC:Results}
The central result of this article is the rate at which the bias of \eqref{EQ:KLKNN_Estimator} decays in the time series context, which we now state.

\begin{theorem}
Fix $k$ in \eqref{EQ:KLKNN_Estimator} and suppose $\{X_i\}_{i=0}^N$ satisfies assumptions (A1), (A2), and (A3) with $\epsilon > \min\{d,1+\sqrt{5}\}$, then there exists a constant $C>0$  and $\theta$ in the interval $$ \big(0,\min \big\{ \frac{\epsilon}{1+\epsilon}, \frac{d+r}{d(2d+r)}, \frac{d+r}{2(2d+r)}, \frac{\epsilon(d+r)}{2(2d+r)(d+1)(2+\epsilon)} \big\}\big) $$ such that 
\begin{equation}
\big| E[\hat{H}_{N+1}] - H(f) \big| \le  \frac{C\log N }{N^\theta}.
\end{equation}
\label{THM:Consistency}
\end{theorem}
There are a few notable implications of this result.  First, the requirement that $\epsilon > \min\{3,1+\sqrt{5}\}$ indicates that the process must mix more quickly in higher dimensions.  Additionally, the mixing condition is trivially satisfied in the \emph{iid} setting so this result extends to that case automatically.  However, even when the stationary density has all moments ($r =\infty$), this result predicts a decay in the bias which is no faster than $O(N^{-1/(d+1)})$ (ignoring the $\log N$ factor) which is slower than the optimal rate from  \cite{delattre_kozachenkoleonenko_2017}.

To prove this result we begin with an integral formulation of the bias following the approach of \cite{delattre_kozachenkoleonenko_2017}. 

\subsection{Integral formulation of the bias}
Note that 
\begin{align*}
&\big| H(f) - E[\hat{H}_{N+1}]\big| = \bigg|E\bigg[\frac{1}{N+1}\sum_{i=0}^N ( -\log f(X_i)  - \log Y_i)\bigg]\bigg| \\
& \le \frac{1}{N+1}\sum_{i=0}^N \bigg| E\big[E[\log T_{k,X_i} -\log Y_i \mid X_i]\big] \bigg| 
\end{align*}
Recall, $T_{k,X_i}\sim \Gamma(k,f(X_i)e^{\Psi(k)})$ so 
\begin{equation}
P(T_{k,X_i} > r ) = \sum_{j=0}^{k-1} e^{-\lambda_{r,k,X_i}}\frac{\lambda_{r,k,X_i}^j}{j!}, \quad \lambda_{r,k,X_i} = e^{\Psi(k)} f(X_i) r.
\label{EQ:Gamma_tail_bound}
\end{equation}  From \cite{delattre_kozachenkoleonenko_2017} Lemma 23, we may write
\begin{equation}
\begin{split}
&E[\log T_{k,X_i} -\log Y_i  \mid X_i] \\
&=\int_0^\infty \big[P(Y_i >r \mid X_i) - P(T_{k,X_i} > r)\big] \frac{dr}{r} \\
&= \int_0^\infty \bigg[ P\bigg(\knn{k}{\mc N_i} > \bigg(\frac{e^{\Psi (k)}  r}{N\nu_d}\bigg)^{1/d} \,\, \bigg\lvert \,\, X_i\bigg)  \\
& \qquad\qquad\qquad\qquad- \sum_{j=0}^{k-1} e^{-\lambda_{r,k,X_i}}\frac{\lambda_{r,k,X_i}^j}{j!}  \bigg] \frac{dr}{r}\\
&= \int_0^\infty \bigg[P(\knn{k}{\mc N_i} > (r/N)^{1/d} \mid X_i) \\
& \qquad\qquad\qquad\qquad- \sum_{j=0}^{k-1}e^{-f(X_i)\nu_d r}\frac{(f(X_i)\nu_d r)^j}{j!} \bigg] \frac{dr}{r}
\end{split}
\label{EQ:Integral_form_cond_bias}
\end{equation}
where the final expression follows from the change of variables $r \mapsto \nu_d e^{-\Psi(k)}r.$ 
For $\xi \in (0,1))$ we will split \eqref{EQ:Integral_form_cond_bias} into the two separate integrals
\begin{align*}
A_N(X_i) &= \int_0^{N^\xi }  \bigg[P(\knn{k}{\mc N_i} > (r/N)^{1/d} \mid X_i)\\
& \qquad \qquad  - \sum_{j=0}^{k-1}e^{-f(X_i)\nu_d r}\frac{(f(X_i)\nu_d r)^j}{j!} \bigg] \frac{dr}{r} \\ 
B_N(X_i) &=  \int_{N^\xi}^\infty \bigg[P(\knn{k}{\mc N_i} > (r/N)^{1/d} \mid X_i) \\
& \qquad \qquad- \sum_{j=0}^{k-1}e^{-f(X_i)\nu_d r}\frac{(f(X_i)\nu_d r)^j}{j!} \bigg] \frac{dr}{r}
\end{align*}
so that we may bound the bias from above by
\begin{equation}
E \big| H(f) - \hat{H}_{N+1}  \big| \le \frac{1}{N+1} \sum_{i=1}^N \big(E|A_N(X_i)| + E |B_N(X_i)|\big).
\end{equation}
We turn our focus to constructing bounds on $E|A_N(X_i)|$ and $E |B_N(X_i)|$ which are uniform in $i$.

\subsection{Bounds on interior terms: $E|A_N(X_i)|$}
\label{SEC:Results.Interior}
Note that $\knn{k}{\mc N_i} > (r/N)^{1/d}$ implies that at most $k-1$ members of the set $\{X_j\}_{j\in \mc N_i}$ are in $B(X_i,(r/N)^{1/d})$.  Additionally, $$\sum_{j=0}^{k-1}e^{-f(X_i)\nu_d r}\frac{(f(X_i)\nu_d r)^j}{j!} = P(Z_{r,X_i} \le k-1)$$ where $Z_{r,X_i} \sim \pois{f(X_i)\nu_d r}.$ Therefore, 
\begin{align*}
A_N(X_i) &= \int_0^{N^\xi}\big[ P(W_{r,X_i,N} \le k-1 \mid X_i)\\
& \qquad\qquad - P(Z_{r,X_i} \le k-1)\big] \frac{dr}{r}
\end{align*}
where 
\begin{equation}
W_{r,X_i,N} = \sum_{j\in \mc N_i} \ind_{B(X_i,(r/N)^{1/d})}(X_j),
\end{equation}
and we may make the upper bound
\begin{equation}
|A_N(X_i)| \le \int_0^{N^\xi} \tv{\ms L(W_{r,X_i,N})\mid X_i), \ms L(Z_{r,X_i })}\frac{dr}{r} .
\label{EQ:A_N_bound}
\end{equation}

By using the total variation distance, \eqref{EQ:A_N_bound} is uniform in $k$.  
To control the total variation between $\ms L (W_{r,X_i,N} \mid X_i)$ and $\ms L(Z_{r,X_i})$ we will make use of the following theorem from \cite{arratia1989}.

\begin{theorem}[\cite{arratia1989}]Let $J$ be an arbitrary index set, and for each $j \in J$, let $\ind_j$ be a Bernoulli random variable with $p_j = P(\ind_j = 1) \in (0,1).$  Let $W = \sum_{j \in J}\ind_j$ be the number of occurrences of dependent events, and let $Z$ be a Poisson random variable with $EZ = EW = \lambda$.  For each $j$, define a neighborhood $B_j \subset J$ containing $j$, then
$$\tv{ \ms L(W), \ms L(Z)} \le 2(b_1+b_2+b_3)$$ where 
\begin{align*}
b_1 &= \sum_{j \in J} \sum_{k\in B_j} p_j p_k,  \\
b_2 &= \sum_{j \in J} \sum_{k \in B_j\setminus \{j\}} P(\ind_j = \ind_k = 1), \\
b_3 &=   \sum_{j\in J}  E\big| E\{\ind_j- p_j\mid \sigma(\ind_k: k \in J - B_j)\} \big|.
\end{align*}
\label{THM:Arratia}
\end{theorem}
In our setting, we'll let $p_j = P(X_j \in B(X_i,(r/N)^{1/d} \mid X_i)$ and use $\mc N_i$ as our index set. Formally, we may think of $X_j$ as depending strongly on $\{X_k\}_{k\in B_j}$ and weakly on $\{X_k\}_{k \notin B_j}$. 

\begin{lemma}
Fix $\beta \in [1/(1+\epsilon), 1)$ and $r > 0$.  Suppose $\{X_i\}_{i=0}^N$ satisfies assumptions (A1) and (A3), then 
\begin{align*}
&\tv{\ms L (W_{r,X_i,N} \mid X_i ), \pois{E W_{r,X_i,N}}}\\
& \le  5(1+K)^2N^{\beta+1}p_{r,X_i,N}^2 + (2L+K) p_{r,X_i,N} 
\end{align*}
where $p_{r,X_i,N} = \int_{B(X_i,(r/N)^{1/d})} f(y) dy$ and $L = \sum_{k=-\infty}^\infty \psi(k).$ 
\label{LEM:TV_Bound}
\end{lemma}

Importantly, for $ 0 \le r \le N^\xi$, $p_{r,X_i,N}$ is shrinking, and we may use Lemma \ref{LEM:TV_Bound} to bound $|A_N(X_i)|.$ 

\begin{lemma} Assume $\{X_i\}_{i=0}^N$ satisfies assumptions (A1) and (A3). For $\xi \in (0,1)$, there exists constants $C_1$, $C_2$, and $C_3$ such that
\begin{equation}
A_N(X_i) \le C_1 N^{2\xi  - \frac{\epsilon}{1+\epsilon} } + C_2N^{\xi-1} +  C_3N^{\frac{(2+d)\xi - 2}{d} }
\label{EQ:Interior_rate}
\end{equation}
so that $E|A_N(X_i)|$ decays whenever $\xi < \frac{2}{2+d} \vee \frac{\epsilon}{2(1+\epsilon)}.$
\label{LEM:Interior_bounds}
\end{lemma}
The $O(N^{\xi - 1})$ term above will always decay fastest so the first and third terms ultimately govern the rate at which $E|A_N(X_i)|$ decays.  Since $2\xi - \frac{1+\epsilon}{2+\epsilon}$ is a decreasing function of $\epsilon$, the first term will decay faster when $\{X_i\}_{i=0}^N$ mixes more rapidly, i.e. $\epsilon$ is larger.  Alternatively, $\frac{(2+d)\xi - 2}{d}$ is bounded from below by $-2/d$ indicating $A_N(X_i)$ decays no faster than $N^{-2/d}$ which is a natural limit when one estimates $Np_{r,X_i,N}$ by $f(X_i)r\nu_d$ under our regularity assumptions on $f$ \cite{delattre_kozachenkoleonenko_2017}.  Unfortunately, the bound on $E|B_N(X_i)|$ place a lower bound on $\xi$ precluding this optimal rate.
Proofs of Lemmas \ref{LEM:TV_Bound} and \ref{LEM:Interior_bounds} are contained in Appendix \ref{APP:Results.Interior.Proof}.

\subsection{Bounds on tail terms: $E|B_N(X_i)|$}
\label{SEC:Results.Tail}
We begin with the trivial bound 
\begin{equation}
\begin{split}
|B_N(X_i)| &\le \int_{N^\xi}^\infty P(\knn{k}{\mc N_i} > (r/N)^{1/d}\mid X_i) \frac{dr}{r}\\
&\qquad + \sum_{j=0}^{k-1}\int_{N^\xi}^\infty e^{-f(X_i)\nu_d r} \frac{(f(X_i)\nu_d r)^j}{j!} \frac{dr}{r}\\
&=\int_{N^\xi}^\infty P(\knn{k}{\mc N_i} > (r/N)^{1/d}\mid X_i ) \frac{dr}{r}\\
& \qquad + \sum_{j=0}^{k-1}\int_{f(X_i) \nu_d N^\xi}^\infty e^{- r} \frac{ r^{j-1} }{j!} dr
\end{split}
\label{EQ:Trivial_tail_bound}
\end{equation}
Note, for any $\theta \in (0,1)$, there exists a constant $C>0$ such that $\sum_{j=0}^{k-1} e^{-r}\frac{r^{j-1}}{j!} \le Cr^{-(1+\theta)}$ for $r > 0$ so that 
$$\sum_{j=0}^{k-1}\int_{f(X_i) \nu_d N^\xi}^\infty e^{- r} \frac{ r^{j-1} }{j!} dr \le \frac{C}{(f(X_i)\nu_d N^\xi)^\theta}.$$
By (A2), choosing $\theta \in \big(0, \frac{d+r}{2d+r}\big)$, it follows that 
$$ \int f(x)  \frac{C}{(f(x)\nu_d N^\xi)^\theta} dx \le \frac{C'}{N^{\theta \xi}}$$
so that 
\begin{equation}
E|B_N(X_i)| \le \frac{C'}{N^{\theta\xi}} + E\bigg[\int_{N^\xi}^\infty P(\knn{k}{\mc N_i} > (r/N)^{1/d}\mid X_i) \frac{dr}{r}\bigg]
\label{EQ:Gamma_tail_knn}
\end{equation} 
Thus, we turn our focus to
\begin{equation}
\int_{N^\xi}^\infty P(\knn{k}{\mc N_i} > (r/N)^{1/d} \mid X_i \big) \frac{dr}{r}
\label{EQ:Tail_knn}
\end{equation}
and present an important bound in terms of the binomial distribution.
\begin{lemma}
\label{LEM:binomial_approx_bound}
Let $p_{r,X_i,N} = \int_{B(X_i,(r/N)^{1/d})} f(y)dy.$ Suppose $\{X_i\}_{i=0}^N$ satisfies assumption (A3). Then for $\beta \in [1/(2+\epsilon),1)$
\begin{equation}
\begin{split}
P(\knn{k}{\mc N_i} &> (r/N)^{1/d} \mid X_i) \\
 & \le e^K\sum_{j=0}^{k-1}  \big(N^{(1-\beta)} p_{r,X_i,N}\big)^j (1- p_{r,X_i,N})^{N^{1-\beta} - j},
\end{split}
\label{EQ:Binomial_bound}
\end{equation}
so that 
\begin{equation}
P(\knn{k}{\mc N_i} > (r/N)^{1/d} \mid X_i) \le C_k N^{k(1-\beta)}e^{-N^{1-\beta}  p_{r,X_i,N}}
\label{EQ:Bound_binomial_bound}
\end{equation} 
where $C_k = ke^{k+K}.$ 
\end{lemma}
The proof of this result is contained in Appendix \ref{APP:Results.Tail.Proof}.   The key observation therein resides in considering the distribution of $\knn{k}{\mc A_i}$ where $\mc A_i$ is some subset $\mc A_i \subset \mc N_i$ so that $\knn{k}{\mc N_i} \le \knn{k}{\mc A_i}$.    By choosing $\mc A_i$ so that all indices are separated by at least $N^\beta$ from $i$ and each other, it follows that $\{X_j\}_{j\in \mc A_i}$ are nearly independent up to a multiplicative correction of the form $1+\psi(N^{\beta}).$

Following the approach of \cite{delattre_kozachenkoleonenko_2017},  let $g(X_i) = 1 + E\| X'- X_i\|^d$ where $X'\sim f$ is independent of $X_i$ and divide \eqref{EQ:Tail_knn} into the separate integrals
\begin{align}
B_{1,1}(X_i) &= \int_{N^\xi}^{2Ng(X_i)} P(Y_i > (r/N)^{1/d} \mid X_i)\frac{dr}{r} \\
B_{1,2}(X_i) &= \int_{2Ng(X_i)}^\infty P(Y_i > (r/N)^{1/d} \mid X_i)\frac{dr}{r} 
\end{align}

\begin{lemma} Suppose $\{X_i\}_{i=0}^N$ satisfies (A2) and (A3), then  for any $\theta >0$, there exists a $C>0$ such that $B_{1,2}(X_i) \le Cg(X_i)N^{-\theta}.$ 
\label{LEM:Tail_bound_knn}
\end{lemma}
This result follows from Lemma \ref{LEM:binomial_approx_bound} and an application of Markov's inequality.  The full details are contained in Appendix \ref{APP:Results.Tail.Proof}.  Importantly, as $g(X_i)$ is integrable by (A2), then for any $\theta >0$, there is $C>0$ such that $E B_{1,2}(X_i) < CN^{-\theta}$ leaving only $B_{1,1}(X_i)$ for consideration. 

As $P( \knn{k}{\mc N_i} >(r/N)^{1/d}\mid X_i)$ is decreasing in $r$, one may attain the trivial bound 
\begin{align*}
B_{1,1}(X_i) &\le P(\knn{k}{\mc N_i} > N^{(\xi-1)/d} \mid X_i ) \int_{N^\xi}^{2Ng(X_i)} \frac{dr}{r} \\
&=P(\knn{k}{\mc N_i} > N^{(\xi-1)/d} \mid X_i ) \\
& \qquad \times \big(\log (2g(X_i)) + (1-\xi)\log(N) \big).
\end{align*}
While $\knn{k}{\mc N_i}$ is a decreasing function of $N$, the rate of decay depends largely on $X_i.$  Formally, we can expect $\knn{k}{\mc N_i}$ to approach zero rapidly when $X_i$ is in a region of high probability ($f(X_i)$ large) but more slowly when $f(X_i)$ is small.  More concretely, when $f(X_i)$ is not too small we may expect $P(\knn{k}{\mc N_i} > N^{(\xi-1)/d} \mid X_i)$ to decay exponentially fast, captured by the following lemma.
\begin{lemma}
In addition to the conditions of Lemma \ref{LEM:binomial_approx_bound}, suppose that $f(X_i) \ge N^{\frac{\epsilon(\xi - 1)}{(d+1)(2+\epsilon)}}$ where $\epsilon$ is taken from (A3), then for $\xi \in (0,1)$
\begin{equation}
\begin{split}
&P(\knn{k}{\mc N_i} > N^{(\xi-1)/d} \mid X_i )\\
& \le C N^{k(1-\beta)}\exp\bigg(- C' N^{1-\beta +\big(1+\frac{\epsilon}{(d+1)(2+\epsilon)}\big)(\xi-1)} \bigg)
\end{split}
\end{equation}
for some constant $C,C'>0$. Furthermore, under (A2) we may conclude that for any $\theta >0$, there is a constant $C''$ such that  $$E [B_{1,1}(X_i) \ind_{f(X_i) \ge N^{\frac{\epsilon(\xi - 1)}{(d+1)(2+\epsilon)}}}(X_i)]\le C''N^{-\theta}$$
whenever $\xi > \frac{d+\epsilon+1}{2(d+\epsilon+1) + d\epsilon}$.
\label{LEM:Middle_bound}
\end{lemma}

Alternatively, we must consider the case where $0< f(X_i) < N^{\frac{\epsilon(\xi - 1)}{(d+1)(2+\epsilon)}}$.

\begin{lemma}
Suppose $0< f(X_i) < N^{\frac{\epsilon(\xi - 1)}{(d+1)(2+\epsilon)}}$.  Assuming (A2), then for any $\theta \in \big(0, (d+r)/(2d+r)\big)$ there exists a constant $C'' >0$ such that  $$E [ B_{1,1}(X_i) \ind_{f(X_i) < N^{(\xi-1)/d)}}(X_i)] \le  N^{\frac{\theta \epsilon(\xi - 1)}{(d+1)(2+\epsilon)}}(\log N + C'').$$
\label{LEM:Mid_bound_decaying_f}
\end{lemma}

In summary, these results imply that $$E|B_N(X_i)| = O\bigg(\max\bigg\{\log (N) N^{\frac{\theta \epsilon(\xi - 1)}{(d+1)(2+\epsilon)}}, N^{-\theta \xi}\bigg\}\bigg)
$$ when $\xi > \frac{d+\epsilon+1}{2(d+\epsilon+1) + d\epsilon}$ and $\theta \in (0,(d+r)/(2d+r))$.  The lower bound on $\xi$ ultimately limits the efficiency of \eqref{EQ:KLKNN_Estimator} one would hope to attain by only examining the contributions from $E|A_N(X_i)|$. 

\subsection{Proof of Theorem \ref{THM:Consistency}}
\begin{proof}
From Lemma \ref{LEM:Interior_bounds},
\begin{equation}
E|A_N(X_i)| = O\bigg( \max \bigg\{N^{2\xi- \frac{\epsilon}{1+\epsilon} },  N^{\frac{(2+d)\xi -2}{d}}\bigg\}\bigg),
\end{equation}
 which is decaying when $\xi < \frac{2}{2+d} \vee \frac{\epsilon}{2(1+\epsilon)}.$  From \eqref{EQ:Gamma_tail_knn} and Lemmas \ref{LEM:Tail_bound_knn}, \ref{LEM:Middle_bound}, and \ref{LEM:Mid_bound_decaying_f},
$$E|B_N(X_i)| = O\bigg( \max\bigg\{\log (N) N^{\frac{\theta\epsilon(\xi-1)}{(d+1)(2+\epsilon)}{d}}, N^{-\theta\xi}\bigg\}\bigg)$$
for $\theta \in (0,(d+r)/(2d+r))$ whenever $\xi > \frac{d+\epsilon+1}{2(d+\epsilon+1) + d\epsilon}.$  

Thus, we require that $$\frac{d+\epsilon+1}{2(d+\epsilon+1) + d\epsilon} < \xi  <  \frac{2}{2+d} \vee \frac{\epsilon}{2(1+\epsilon)}$$
 which holds whenever $\epsilon >  \min \{d ,  1+\sqrt{5}\}.$
 
Importantly, sending $\epsilon, r \to \infty$ one can see that $\xi \in (0,1/2)$. Inserting these extremal values of $\xi$ into the bounds on $E|A_N(X_i)|$ and $E|B_N(X_i)|$ gives the interval from Thm. \ref{THM:Consistency}. 
\end{proof}

\section{Example: Stationary Gaussian Distribution}
\label{SEC:Example}
Let $0_d$ denote the zero vector in $\mathbb{R}^d$, $I_d$ denote the identity matrix in $\mathbb{R}^{d\times d}$, and $$\Sigma_d = \begin{bmatrix}
1 & r & 0 &\dots & 0 \\ 
r  & 1 & r  & \ddots &  \vdots\\
0 & \ddots & \ddots & \ddots & 0 \\
\vdots & \ddots & \ddots & \ddots & r  \\
0 &\dots &0  & r & 1  
\end{bmatrix}$$
for $r \in (-1,1)$. Let $\mc N(\mu,\Sigma)$ denote the (multivariate) Gaussian distribution with mean $\mu$ and (co)variance $\Sigma.$
Suppose that $[X_{i+1},X_i]$ has joint distribution
\begin{equation}
\begin{bmatrix}
X_{i+1} \\ X_i
\end{bmatrix} 
\sim   N\left(
\begin{bmatrix}
0_d \\ 
0_d
\end{bmatrix}, 
\begin{bmatrix}
\Sigma_d & \rho I_d \\ 
\rho I_d& \Sigma_d
\end{bmatrix}\right)
\label{EQ:Stationary_Joint}
\end{equation}
for $\rho \in (-1,1)$  
so that $\{X_i\}_{i=1}^N$ is a Markov chain with stationary distribution $\mc N (\mu_d,\Sigma_d)$. The (transition) conditional distribution of $X_{i+1}|X_i$ is
\begin{equation}
X_{i+1}|X_i \sim  \mc N \left( \rho \Sigma_d^{-1} 
 X_i , \Sigma_d - \rho^2 \Sigma_d^{-1}  \right).
 \label{EQ:Gaussian_MC}
 \end{equation}
The true entropy of the stationary Gaussian distribution is
\begin{equation}
H(X_1) = \frac{d}{2}+\frac{d}{2}\log 2\pi + \frac{1}{2}\log |\Sigma_d|
\label{EQ:Entropy_Normal}
\end{equation}
where the determinant of $\Sigma_d$ satisfies the recurrence relation $$|\Sigma_d| = |\Sigma_{d-1}| - r^2 |\Sigma_{d-2}|, \quad |\Sigma_1| = 1, |\Sigma_2| = 1-r^2.$$

\begin{figure}
\begin{center}
\adjincludegraphics[ trim={0 0 {.5\width} 0}, clip, width = 3 in]{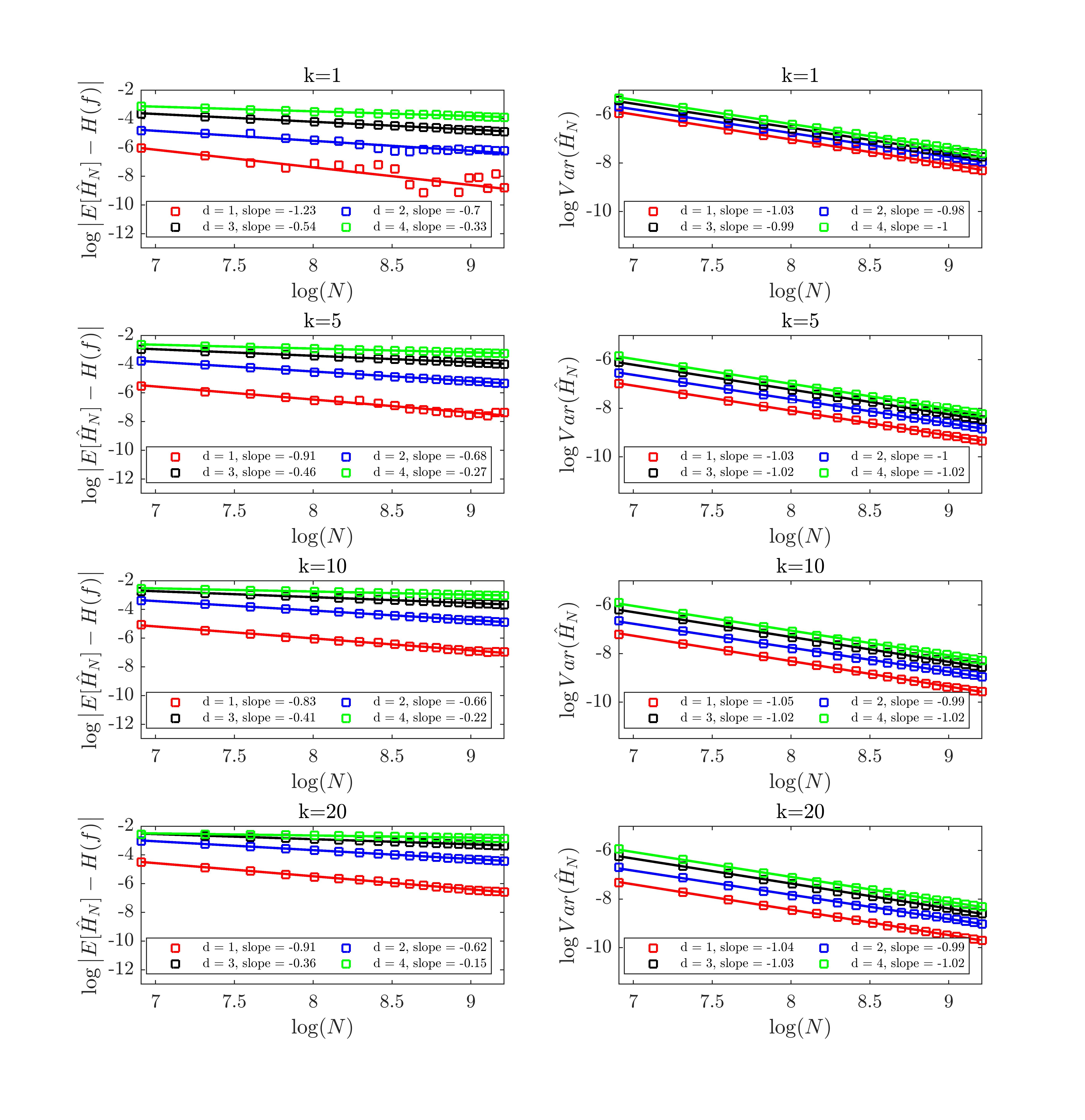}
\end{center}
\caption{Numerical estimates of the bias  of \eqref{EQ:KLKNN_Estimator} applied to \eqref{EQ:Gaussian_MC} with entropy \eqref{EQ:Entropy_Normal} as a function of the sample size are shown above.  In each case, 2000 iterations of the Markov Chain were generated taking $r=\rho = 1/4$ to estimate $E[\hat{H}]$ and $Var(\hat{H}_N)$. In each frame, the squares correspond to the logarithm of the estimated bias  and the lines correspond to least-squares fits in log-log space.  The slopes of these lines, rounded to two decimal places, are included in the legend for reference to identify the rate at which these quantities decay.}
\label{FIG:Bias_r25_p25}
\end{figure}

To study the efficiency of the $k$-nearest neighbor entropy estimator, we generated 2000 realizations of this stationary process with $r = 1/4$ and $\rho = 1/4$ and compared the mean estimate with the true entropy for various sample lengths $N$, dimensions $d$, and number of nearest neighbors $k$.  The Euclidean norm was used in all distance calculations.  The results for the bias are shown in Fig. \ref{FIG:Bias_r25_p25} including variance estimates in Fig. \ref{FIG:Variance_r25_p25}.  

\begin{figure}
\begin{center}
\adjincludegraphics[ trim={{.5\width} 0 0  0}, clip, width = 3 in]{Bias_variance_Gauss_MC_r25e-1_p25e-1.jpg}
\end{center}
\caption{Numerical estimates of the variance of \eqref{EQ:KLKNN_Estimator} applied to \eqref{EQ:Gaussian_MC} with entropy \eqref{EQ:Entropy_Normal} as a function of the sample size are shown above.  In each case, 2000 iterations of the Markov Chain were generated taking $r=\rho = 1/4$ to estimate $E[\hat{H}]$ and $Var(\hat{H}_N)$. In each frame, the squares correspond to the logarithm of the estimated  variance and the lines correspond to least-squares fits in log-log space.  The slopes of these lines, rounded to two decimal places, are included in the legend for reference to identify the rate at which these quantities decay.}
\label{FIG:Variance_r25_p25}
\end{figure}

The stationary Gaussian density is Lipschitz and has moments of all orders.  Additionally, the autoregressive nature of the process with $ - 1<\rho  <1$ can be expected to mix exponentially fast.  Thus, we may formally treat $\alpha=1$ in (A1) and $r$ and $\epsilon$ as infinite from assumptions (A2) and (A3).  In this case, by Thm. \ref{THM:Consistency} the maximal rate of decay of the bias is $O(N^{-1/(d+1)})$ up to $\log N$.  This is consistent with the simulation results (Figure \ref{FIG:Bias_r25_p25}) which all show the bias decaying at faster rates than predicted.    While not considered analytically, the numerical results show the variance of $\hat{H}_N$ is are approximately $O(N^{-1})$ for all dimensions (Figure \ref{FIG:Variance_r25_p25}).  Both \cite{berrett2019} and \cite{delattre_kozachenkoleonenko_2017} identified a distribution-free inflation of the variance which one should expect to extend to the time series setting. However, the numerical results were too noisy to verify this fact.

\appendix
\section{Proofs of Lemmas}
\subsection{Proofs for Section \ref{SEC:Results.Interior}}
\label{APP:Results.Interior.Proof}
\begin{proof} (Lemma \ref{LEM:TV_Bound})
Note that (A3) implies 
\begin{equation}
\begin{split}
&P(X_j \in B(X_i,(r/N)^{1/d} \mid X_i )\\
& \le [1+\psi( |i-j|) ] p_{r,X_i,N} \le [1+K] p_{r,X_i,N}
\end{split}
\label{EQ:Mixing_bound_v2}
\end{equation}
Let $B_j = \mc N_i \cap\{j - \floor{N^\beta},\dots, j+\floor{N^\beta} \}$ so that $|B_j| \le  2N^\beta+1$.  Then 
\begin{align*}
b_1 & = \sum_{j\in \mc N_i} \sum_{k\in B_j} P(X_j \in B(X_i,(r/N)^{1/d}) \mid X_i)\\
& \qquad  \qquad  \qquad \times P(X_k \in B(X_i,(r/N)^{1/d})\mid X_i) \\
&\le \sum_{j\in \mc N_i}\sum_{k\in B_j}[1+K]^2p_{r,X_i,N}^2\\
& = (2N^\beta+1)(1+K)^2Np_{r,X_i,N}^2\\
& \le 3(1+K)^2N^{1+\beta}p_{r,X_i,N}^2
\label{EQ:b1_bound}
\end{align*} 
For $b_2$, we again use (A3) to attain the bound
\begin{align*}
& P\big(X_k \in B(X_i,(r/N)^{1/d}) \mid X_i, X_j \in B(X_i,(r/N)^{1/d})\big)  \\
&\le [1+\psi( \min \{|i-j|,|j-k|\})] p_{r,X_i,N} \le [1+K]p_{r,X_i,N}
\end{align*}
\begin{align*}
b_2 &= \sum_{j\in \mc N_i} \sum_{k\in B_j\setminus\{j\}} P\big(X_j \in B(X_i,(r/N)^{1/d}) \\
& \qquad \qquad \qquad \qquad \times P( X_k \in B(X_i,(r/N)^{1/d}) \mid X_i \big) \\
&= \sum_{j\in \mc N_i} \sum_{k\in B_j\setminus\{j\} } [1+K]^2p_{r,X_i,N}^2 =2N^{\beta+1} (1+K)^2 p_{r,X_i,N}^2
\end{align*}
Finally, 
\begin{align*}
& E\bigg| E\bigg[\ind_{B(X_i,(r/N)^{1/d})}(X_j) \\
& \qquad \qquad \qquad- P\big(X_j \in B(X_i,(r/N)^{1/d} \mid X_i\big) \,\bigg|\,\mc F_{\mc N_i \setminus B_j} \bigg]\bigg|\\
&\le \bigg| P\big(X_j \in B(X_i,(r/N)^{1/d}) \mid X_i \big) \\
& \qquad\qquad- P\big(X_j \in B(X_i,(r/N)^{1/d})\big)  \mid X_i, \mc F_{\mc N_i \setminus B_j }) \bigg| \\
 &\le \big[\psi(| i- j |) + \psi\big(d(j, \{0,\dots, N\} \setminus B_j)\big)\big] p_{r,X_i,N} \\
 &\le \big(2\psi(|i-j|) + \psi(\floor{N^\beta})\big) p_{r,X_i,N}.
\end{align*}
Using these results,
\begin{align*}
b_3 &= \sum_{j \in \mc N_i} E\bigg| E\big[\ind_{B(X_i,(r/N)^{1/d})}(X_j) \\
& \qquad\qquad- P(X_j \in B(X_i,(r/N)^{1/d} \mid X_i) \mid \mc F_{\mc N_i \setminus B_j} \big]\bigg| \\
& \le p_{r,X_i,N} \sum_{j\in \mc N_i} \big(2\psi(|i-j|) + \psi(\floor{N^\beta})\big)
\end{align*}
From (A3), $L =\sum_{k=-\infty}^\infty \psi(k) < \infty$, and 
$
b_3 \le \bigg(2L + \frac{KN}{1+N^{\beta(1+\epsilon)}}\bigg) p_{r,X_i,N}.
$ Since $\beta \ge 1/(1+\epsilon)$, 
$$\frac{N}{1+N^{\beta(1+\epsilon)}} \le N^{1-\beta(1+\epsilon)} \le 1$$ so $b_3 \le (2L+K)p_{r,X_i,N}.$
\end{proof}

\begin{proof} (Lemma \ref{LEM:Interior_bounds})
Throughout this proof we set $\beta = 1/(1+\epsilon).$ Now $|A_N(X_i)|$ is bounded above by 
\begin{equation}
\begin{split}
&  \int_0^{N^\xi} \tv{\ms L(W_{r,X_i,N}\mid X_i),\pois{EW_{r,X_i,N} } } \frac{dr}{r} \\
&  + \int_0^{N^\xi} \tv{\pois{EW_{r,X_i,N} }, \pois{f(X_i)\nu_d r} }  \frac{dr}{r}
\label{EQ:TV_bound_split}
\end{split}
\end{equation} 
From Lemma \ref{LEM:TV_Bound},
\begin{align*}
&\int_0^{N^\xi} \tv{\ms L(W_{r,X_i,N}\mid X_i),\pois{EW_{r,X_i,N} } } \frac{dr}{r}\\
  &\le \int_0^{N^\xi} \big(5(1+K)^2N^{\beta+1}p_{r,X_i,N}^2 + (2L+K) p_{r,X_i,N} \big) \frac{dr}{r} 
\end{align*}
From (A1) and \eqref{EQ:Holder_bound_on_integral} 
\begin{align*}
p_{r,X_i,N} &\le f(X_i)\nu_d \frac{r}{N} + C_f\nu_d \bigg(\frac{r}{N}\bigg)^{(2+d)/d}\\
 &= \frac{\nu_d r}{N}\bigg(f(X_i) + C_f \bigg(\frac{r}{N}\bigg)^{2/d}\bigg) \le  \frac{2C_f\nu_d r}{N} 
 \end{align*}
 for $0\le r \le N^{\xi}$. Thus, 
\begin{equation}
\begin{split}
&\int_0^{N^\xi} \big(5(1+K)^2N^{\beta+1}p_{r,X_i,N}^2 + \underbar{C} p_{r,X_i,N} \big) \frac{dr}{r}  \\
& \le \frac{20(1+K)^2\nu_d^2C_f^2}{N^{1-\beta}} \int_0^{N^\xi} r dr + \frac{2C_f (2L+K)\nu_d}{N} \int_0^{N^\xi} dr \\
&= C_1 N^{2\xi +\beta -1 } + C_2 N^{\xi - 1} = C_1N^{2\xi - \frac{\epsilon}{1+\epsilon}} + C_2N^{\xi - 1}
\end{split}\label{EQ:Bound11}
\end{equation}
where $C_1$ and $C_2$ are bounded constants which do not depend on $X_i$. 

Note that $\tv{\pois{\mu_1},\pois{\mu_2}} \le |\mu_1-\mu_2|.$ Thus,
\begin{align*}
& \tv{\pois{E[W_{r,X_i,N}\mid X_i]  }, \pois{f(X_i)\nu_d r} }  \\
&\le  \bigg| \sum_{j \in \mc N_i} P\big( X_j \in B(X_i,(r/N)^{1/d}) \mid X_i \big) - f(X_i) \nu_d r  \bigg| \\
&\le \sum_{j \in \mc N_i} \big| P\big( X_j \in B(X_i,(r/N)^{1/d}) \mid X_i \big) - p_{r,X_i,N} \big| \\
& \qquad\qquad\qquad+ |Np_{r,X_i,N} - f(X_i)\nu_d r| 
\end{align*}
Making use of (A3), it follows that 
\begin{align*}
& \sum_{j \in \mc N_i} \big| P\big( X_j \in B(X_i,(r/N)^{1/d}) \mid X_i \big) - p_{r,X_i,N} \big|\\
& \le \sum_{k=-\infty}^\infty \psi(k) p_{r,X_i,N} = Lp_{r,X_i,N}.
\end{align*}
Again from (A1) and \eqref{EQ:Holder_bound_on_integral},
$$ |Np_{r,X_i,N} - f(X_i)\nu_d r| \le C_f\nu_d \frac{r^{(2+d)/d}}{N^{2/d}}.$$
Applying these results to the last term in \eqref{EQ:TV_bound_split}, and using the bound on $p_{r,X_i,N}$, 
\begin{equation}
\begin{split}
&\int_0^{N^\xi} \tv{\pois{EW_{r,X_i,N} }, \pois{f(X_i)\nu_d r} }  \frac{dr}{r} \\
&\le \frac{L\nu_d(f(X_i+\overline{C})}{N}\int_0^{N^\xi}dr + \frac{C_f\nu_d}{N^{2/d}} \int_0^{N^\xi} r^{2/d} dr\\
& \le C_3 N^{\xi - 1} + C_4 N^{\frac{2+d}{d}\xi - \frac{2}{d} } 
\end{split}
\label{EQ:Bound12}
\end{equation}
where $C_3$ and $C_4$ are bounded constants which do not depend on $X_i.$   Combining the $O(N^{\xi-1})$ terms in \eqref{EQ:Bound11} and \eqref{EQ:Bound12} yields \eqref{EQ:Interior_rate}.  For this bound to decay in $N$ we then need $\xi  < \frac{\epsilon}{2(1+\epsilon)}$, $\xi < 1$, and $ \xi < \frac{2}{2+d}$ which is equivalent to the condition $\xi < \frac{2}{2+d} \vee \frac{\epsilon}{2(1+\epsilon)}$. 
\end{proof}

\subsection{Proofs for Section \ref{SEC:Results.Tail}}
\label{APP:Results.Tail.Proof}
\begin{proof} (Lemma \ref{LEM:binomial_approx_bound})
Given any subset $\mc A_i\subset \mc N_i$ of indices,  $\rho_{i,X_i,A_i} \ge \rho_{i,X_i,N_i} $
so that 
$$P(\knn{k}{\mc N_i} >r \mid X_i ) \le P(\knn{k}{\mc A_i} > r \mid X_i).$$ 
Take $$\mc A_i = \mc N_i \cap \{\dots, i-2\ceil{N^{\beta}}, i - \ceil{N^{\beta}},i+\ceil{N^{\beta}}, i+2\ceil{N^{\beta}},\dots\}$$
so that $|i_1-i_2| > N^\beta$ for $i_1,i_2\in \mc A_i$ and 
\begin{equation}
\begin{split}
&P(X_{i_1} \in B(X_i,r) \mid X_i, \mc F_{A_{i}\setminus \{i_1 \} } )\\
 &\qquad\qquad\le (1+\psi(N^\beta))P(X_{i_1}  \in B(X_i,r) ), \\
& P(X_{i_1} \notin B(X_i,r) \mid X_i, \mc F_{A_{i}\setminus \{i_1 \} } )  \\
& \qquad\qquad\le (1+\psi(N^\beta))P(X_{i_1}  \notin B(X_i,r) ).
\end{split} 
\label{EQ:Mixing_prob_bounds}
\end{equation}
Thus, $\sum_{j\in\mc A_i} \ind_{B(X_,(r/N)^{1/d})}(X_j)$ follows a $Bin(|\mc A_i|,p_{r,X_i,N})$ distribution up to a multiplicative correction of the form $(1+\psi(N^\beta))^{|\mc A_i|}.$
More explicitly,
\begin{equation}
\begin{split}
&P\big(\knn{k}{\mc A_i} > (r/N)^{1/d} \mid X_i\big) \\
& \le [1+\psi(N^\beta)]^{|\mc A_i|} \sum_{j=0}^{k-1}\binom{|A_i|}{j} p_{r,X_i,N}^j (1-p_{r,X_i,N})^{|\mc A_i| -j}.
\end{split}
\label{EQ:Almost_binom_bound}
\end{equation}
Additionally, as $|\mc A_i| \le N^{1-\beta}$ and $1-\beta(2+\epsilon) < 0$ by assumption
\begin{equation}
\begin{split}
[1+\psi(N^\beta)]^{|\mc A_i|}& \le e^{N^{1-\beta} \log(1+ \psi(N^\beta)) }\\
& \le \exp \bigg(\frac{KN^{1-\beta}}{1+N^{\beta(1+\epsilon)}} \bigg) \le e^K
\end{split}
\label{EQ:Mixing_Rate_Saturation_Bound}
\end{equation}
and $\binom{|\mc A_i|}{j} \le |\mc A_i|^j \le N^{j(1-\beta)}$. Applying these inequalities to \eqref{EQ:Almost_binom_bound} gives \eqref{EQ:Binomial_bound}.  Finally, bounding $p_{r,X_i,N}$ above by one, we have 
\begin{align*}
&p_{r,X_i,N}^j(1-p_{r,X_i,N})^{|\mc A_i| - j}\\
&  \le \exp\big( (N^{1-\beta}-j) \log (1-p_{r,X_i,N}) \big)\\
&\le \exp\big(-(N^{1-\beta}-j) p_{r,X_i,N}\big) \le \exp\big(-N^{1-\beta}p_{r,X_i,N} + k\big)
\end{align*}
where we have used the requirement that $j \le k$.  This bound holds for all $j = 0, \dots, k-1$. Applying it to each term in the sum from \eqref{EQ:Almost_binom_bound} yields the desired result.  
\end{proof}

\begin{proof} (Lemma \ref{LEM:Tail_bound_knn})
We begin with the change of variables $r \mapsto rN$, so that 
$$B_{1,2}(X_i) = \int_{2g(X_i)}^\infty P(\knn{k}{\mc N_i} > r^{1/d} \mid X_i) \frac{dr}{r}.$$
Taking $\beta = 1/(2+\epsilon)$ and applying Lemma \ref{LEM:binomial_approx_bound} we have
\begin{equation}
\begin{split}
B_{1,2}(X_i)  &\le e^K \sum_{j=0}^{k-1} N^{j(1-\beta)}  \\
& \qquad \times \int_{2g(X_i)}^\infty P(X' \notin B(X_i,r^{1/d}) ) ^{\mc A_i - j} \frac{dr}{r}.
\end{split}
\end{equation}
Here, $X'\in \mathbb{R}^d$ is a random variable with density $f$ independent of $X_i.$
From Markov's inequality, $$P(X' \notin B(X_i, r^{1/d})) = P(\| X'-X_i\|^d > r) \le \frac{E\|X'-X_i\|^d}{r}.$$
Then,
\begin{align*}
&B_{1,2}(X_i) \le e^K\sum_{j=0}^{k-1}N^{j(1-\beta)} g(X_i) \\
& \qquad\qquad  \qquad\times \int_{2g(X_i)}^\infty \bigg(\frac{E_{X'\sim f}\|X'-X_i\|}{r}\bigg)^{|\mc A_i|-j-1} \frac{dr}{r^2}\\
&\le e^{K}g(X_i) \sum_{j=0}^{k-1} N^{j(1-\beta)}\bigg(\frac{1}{2}\bigg)^{|\mc A_i|-j-1} 
\end{align*}
The final bound follows as $E_{X'\sim f}\|X'-X_i\|/r \le 1/2$ for $r > 2g(X_i).$ Finally, for any $\theta >0$, there exists a constant $C>0$ such that $$\sum_{j=0}^{k-1}N^{j(1-\beta)}\bigg(\frac{1}{2}\bigg)^{N^{1-\beta}-j-1} \le CN^{-\theta}$$
and we may conclude that $B_{1,2}(X_i) \le Cg(X_i) N^{-\theta}.$ 
\end{proof}

\begin{proof}(Lemma \ref{LEM:Middle_bound}) From Lemma \ref{LEM:binomial_approx_bound}, it follows that 
\begin{align*}
& P(\knn{k}{\mc N_i} > N^{(\xi-1)/d} \mid X_i )\\
& \qquad \le C_k N^{k(1-\beta)}\exp\big\{ - N^{1-\beta}   p_{N^\xi,X_i,X} \big\}.
\end{align*}
The above expression decays exponentially if $N^{1-\beta}p_{N^\xi,X_i,N} > N^{\epsilon'}$ for some $\epsilon'>0$.  We may apply \eqref{EQ:Holder_bound_on_integral} to attain a lower bound on $p_{N^\xi,X_i,N}$,
\begin{align}
& p_{N^\xi,X_i,N} \ge f(X_i)\nu_d N^{\xi-1}\\
&\qquad\qquad\qquad - \frac{d^2\nu_dC_f}{2} N^{(\xi-1)\frac{2+d}{d}}  - C_f\nu_d N^{(\xi-1)\frac{2+d+\alpha}{d}} \\
&= \nu_d N^{\xi-1} \bigg(f(X_i) - \frac{d^2 C_f}{2}N^{\frac{2}{d}(\xi-1)} - C_f N^{(\xi-1)\frac{2+\alpha}{d}}  \bigg). \label{EQ:Lower_bound_extreme}
\end{align}
Assuming $f(X_i) \ge N^{\frac{\epsilon(\xi - 1)}{(d+1)(2+\epsilon)}}$ then
\begin{align*}
& f(X_i) - \frac{d^2 C_f}{2}N^{\frac{2}{d}(\xi-1)} - C_f N^{(\xi-1)\frac{2+\alpha}{d}}\\
&\ge N^{\frac{\epsilon(\xi - 1)}{(d+1)(2+\epsilon)}}  \bigg(1 - \frac{d^2C_f}{2} N^{(\xi-1)(\frac{2}{d} - \frac{\epsilon}{(d+1)(2+\epsilon)} ) } \\
& \qquad\qquad\qquad\qquad- C_f N^{(\xi-1)(\frac{2+\alpha}{d} - \frac{\epsilon}{(d+1)(2+\epsilon)}
)} \bigg)
\end{align*}
Importantly, $\frac{\epsilon}{(d+1)(2+\epsilon)} < \frac{1}{d+1} < \frac{1}{d}$. We may choose
$N$ sufficiently large so that the higher order terms within the parentheses of the last line of the preceding expression are both less than $\frac{1}{3}.$  This gives the lower bound $p_{N^\xi,X_i,N} \ge \frac{1}{3}\nu_d N^{(\xi-1)\big(1+\frac{\epsilon}{(d+1)(2+\epsilon) }\big)}.$  Inserting this quantity into \eqref{EQ:Binomial_bound}, we have
\begin{align*}
&P(\knn{k}{\mc N_i} > N^{(\xi-1)/d} \mid X_i)\\
& \le C_kN^{k(1-\beta)}\exp\bigg(-\frac{1}{3}\nu_d N^{1-\beta + (\xi-1)\big(1+\frac{\epsilon}{(d+1)(2+\epsilon) }\big)} \bigg).
\end{align*}

Now assume $1-\beta + (\xi-1)\big(1+\frac{\epsilon}{(d+1)(2+\epsilon) }\big)> 0$, or equivalently, $\xi > \frac{d+\epsilon+1}{2(d+\epsilon+1) + d\epsilon}$. In this case, for any $\theta >0$, there is a constant $C'>0$ such that the above expression is less than $$P(\knn{k}{\mc N_i} > N^{(\xi-1)/d} \mid X_i ) \le \frac{C'}{N^\theta \log N}.$$ As a result, 
\begin{align*}
&E[B_{1,1}(X_i) \ind_{f(X_i) \ge N^{\frac{\xi-1}{d}}}(X_i)]\\
 &\le \frac{C'}{N^\theta \log N}\int \log \big(2g(X_i) N^{1-\xi} \big) f(X_i) dX_i \\
&\le \frac{C'E[1+2g(X_i)]}{N^\theta \log N} + \frac{C'(1-\xi)}{N^\theta} \le \frac{C''}{N^{\theta}}
\end{align*}
for some $C''>0$ since $\log 2g(X_i) \le 1+2\|X'-X_i\|^d$ which is integrable by (A2).
\end{proof}

\begin{proof}(Lemma \ref{LEM:Mid_bound_decaying_f}) Assuming $0< f(X_i) <N^{\frac{\epsilon(\xi - 1)}{(d+1)(2+\epsilon)}}$ it follows that $$P(\knn{k}{\mc N_i} > N^{(\xi-1)/d} \mid X_i ) \le 1 < \frac{N^{\frac{\theta \epsilon(\xi - 1)}{(d+1)(2+\epsilon)}}}{f^\theta(X_i)}$$
for any $\theta >0.$ Thus, 
\begin{align*}
&E[B_{1,1}(X_i)\ind_{f(X_i) < N^{\frac{\theta \epsilon(\xi - 1)}{(d+1)(2+\epsilon)}}}(X_i) ]\\
& \le  N^{\frac{\epsilon(\xi - 1)}{(d+1)(2+\epsilon)}}\int\log \big(2g(X_i) N^{1-\xi} \big)f^{1-\theta}(X_i)dX_i \\
&\le  \log (N)N^{\frac{\theta \epsilon(\xi - 1)}{(d+1)(2+\epsilon)}} \\
& \qquad\qquad+ N^{\frac{\epsilon(\xi - 1)}{(d+1)(2+\epsilon)}}\int \log (2g(X_i)) f^{1-\theta}(X_i)dX_i.
\end{align*}
We must choose $\theta$ sufficiently small so that the final integral in  the preceding expression is finite.  Consider $r$ from (A2). Let $M >1$ and  set $\theta = \frac{(M-1)(d+r)}{M(2d+r)} < 1$. Then, there exists a constant $C>0$ depending on $M$ such that $\log 2g(X_i) \le C^{1-\theta}(1+\|X_i\|^d)^{\frac{r(1-\theta)}{dM}}$.  As such,
\begin{align*}
&\int \log (2g(X_i)) f^{1-\theta}(x)dx \\
& \le \int  \frac{\big(C(1+\|x\|^d)^{\frac{r}{dM}}(1+\|x\|)^{\frac{Md+(M-1)r}{M} }f(x)\big)^{1-\theta}}{(1+\|x \|)^{(1-\theta)\frac{Md+(M-1)r}{M}}}dx \\
&\le \bigg(C\int (1+\|x\|^d)^{\frac{r}{Md}}(1+\|x\|)^{\frac{Md+(M-1)r}{M} }f(x) dx  \bigg)^{1-\theta} \\
& \qquad \qquad\times \bigg(\int \frac{1}{(1+\|x\| )^{\frac{(1-\theta)(Md+(M-1)r)}{M\theta} }} dx \bigg)^{\theta}.
\end{align*}
by H\"older's inequality.  Since $(1+\|x\|^d)^{\frac{r}{Md}}(1+\|x\|)^{\frac{Md+(M-1)r}{M} }$ is $O(\|x\|^{r+d})$ as $\|x\| \to \infty,$ the first term above is integrable by (A2).  Secondly, recalling $\theta = \frac{(M-1)(d+r)}{M(2d+r)}$, it follows that $(1-\theta)(Md+(M-1)r)/(2\theta) > \frac{M+1}{M}d + \frac{r}{M} > d$ so the second integral is finite.  This holds for any $M >1$.  As a result, for any $\theta \in (0, (d+r)/(2d+r))$ there exists as constant $C''> 0$ such that 
$$E[B_{1,1}(X_i)\ind_{f(X_i) < N^{\frac{\xi-1}{d}}}(X_i) ]  \le N^{\frac{\theta \epsilon(\xi - 1)}{(d+1)(2+\epsilon)}}(\log N + C'').$$
\end{proof}

\section*{Acknowledgements}
The authors would like to thank Jonathan Mattingly for his insight on a number of technical details regarding the tail bound on $P(\knn{k}{\mc N_i} > r \mid X_i)$. ALY was supported by grants \#1045153 and  \#1546130 from the National Science Foundation.  



%\bibliographystyle{IEEEtran}
\bibliographystyle{IEEEtran}
\bibliography{IEEEabrv,references}

\end{document}